# Smoothed weighted empirical likelihood ratio confidence intervals for quantiles

JIAN-JIAN REN

*Department of Mathematics, University of Central Florida, Orlando, Florida 32816, USA.*
*E-mail: jren@mail.ucf.edu*

Thus far, likelihood-based interval estimates for quantiles have not been studied in the literature on interval censored case 2 data and partly interval censored data, and, in this context, the use of smoothing has not been considered for any type of censored data. This article constructs smoothed weighted empirical likelihood ratio confidence intervals (WELRCI) for quantiles in a unified framework for various types of censored data, including right censored data, doubly censored data, interval censored data and partly interval censored data. The fourth order expansion of the weighted empirical log-likelihood ratio is derived and the *theoretical* coverage accuracy equation for the proposed WELRCI is established, which generally guarantees at least '*first order*' accuracy. In particular, for right censored data, we show that the coverage accuracy is at least $O(n^{-1/2})$ and our simulation studies show that in comparison with empirical likelihood-based methods, the smoothing used in WELRCI generally provides a shorter confidence interval with comparable coverage accuracy. For interval censored data, it is interesting to find that with an adjusted rate $n^{-1/3}$, the weighted empirical log-likelihood ratio has an asymptotic distribution completely different from that obtained by the empirical likelihood approach and the resulting WELRCI perform favorably in the available comparison simulation studies.

## 1. Introduction

Since Owen (1988), the empirical likelihood method has been developed to construct tests and confidence sets based on the nonparametric likelihood ratio; see Owen (1990, 1991, 2001), DiCiccio, Hall and Romano (1991), Qin and Lawless (1994), Mykland (1995) and Zhou (2005), among others. Studies have shown that the empirical log-likelihood ratio usually has an asymptotic chi-squared distribution and that the empirical likelihood ratio inference is of comparable accuracy to alternative methods. In particular, it is shown that the empirical likelihood is Bartlett-correctable for smooth function models (DiCiccio, Hall and Romano (1991)).

In survival analysis, the quantiles of a lifetime distribution are often of significant interest. It is known that the likelihood-based methods perform favorably in catching the skewness of the distribution of the statistics of interest (Owen (1988), Ren (2001)) and Chen and Hall (1993) showed that for the complete sample case, smoothing can improve the coverage accuracy of empirical likelihood-based confidence intervals for quantiles. But, thus far, likelihood-based interval estimates for quantiles have not been studied in literature for interval censored case 2 data and partly interval censored data, and, in this context, the use of smoothing has not been considered for any type of censored data. This article studies one type of smoothed weighted empirical likelihood-based interval estimate for the $q$th quantile of the lifetime distribution function $F_0$:

$$\theta_0 = F_0^{-1}(q), \qquad 0 < q < 1, \tag{1.1}$$

with various types of censored data. Throughout this paper, we let $X_1, \ldots, X_n$ be an independently and identically distributed (i.i.d.) random sample from a continuous and non-negative distribution function (d.f.) $F_0$, but we consider the cases when such an i.i.d. sample is not completely observable due to censoring. Specifically, in this work, we consider the following types of censored data.

***Right censored sample.*** The observed data are $\boldsymbol{O}_i = (V_i, \delta_i), i = 1, \ldots, n$, with

$$V_i = \begin{cases} X_i, & \text{if } X_i \le Y_i, \quad \delta_i = 1, \\ Y_i, & \text{if } X_i > Y_i, \quad \delta_i = 0, \end{cases} \tag{1.2}$$

where $Y_i$ is the right censoring variable and is independent of $X_i$. This type of censoring has been extensively studied in the literature in the past few decades.

***Doubly censored sample.*** The observed data are $\boldsymbol{O}_i = (V_i, \delta_i), i = 1, \ldots, n$, with

$$V_i = \begin{cases} X_i, & \text{if } Z_i < X_i \le Y_i, & \delta_i = 1, \\ Y_i, & \text{if } X_i > Y_i, & \delta_i = 2, \\ Z_i, & \text{if } X_i \le Z_i, & \delta_i = 3, \end{cases} \tag{1.3}$$

where $Y_i$ and $Z_i$ are right and left censoring variables, respectively, and are independent of $X_i$ with $P\{Z_i < Y_i\} = 1$. This type of censoring has been considered by Turnbull (1974), Chang and Yang (1987), Gu and Zhang (1993), Ren (1995) and Mykland and Ren (1996), among others. One recent example of doubly censored data was encountered in a study of primary breast cancer (Ren and Peer (2000)).

***Interval censored sample.***

***Case 1.*** The observed data are $\boldsymbol{O}_i = (Y_i, \delta_i), i = 1, \ldots, n$, with

$$\delta_i = I\{X_i \le Y_i\}. \tag{1.4}$$

***Case 2.*** The observed data are $\boldsymbol{O}_i = (Y_i, Z_i, \delta_i), i = 1, \ldots, n$, with

$$\delta_i = \begin{cases} 1, & \text{if } Z_i < X_i \leq Y_i, \\ 2, & \text{if } X_i > Y_i, \\ 3, & \text{if } X_i \leq Z_i, \end{cases} \quad (1.5)$$

where $Y_i$ and $Z_i$ are independent of $X_i$ and satisfy $P\{Z_i < Y_i\} = 1$ for Case 2. These two types of interval censoring were considered by Groeneboom and Wellner (1992), among others. In practice, interval censored Case 2 data were encountered in AIDS research (Kim, De Gruttola and Lagakos (1993); also see the discussion in Ren (2003)).

***Partly interval censored sample.***

***'Case 1' partly interval censored data.*** The observed data are

$$\boldsymbol{O}_i = \begin{cases} X_i, & \text{if } 1 \leq i \leq n_1, \\ (Y_i, \delta_i), & \text{if } n_1 + 1 \leq i \leq n, \end{cases} \quad (1.6)$$

where $\delta_i = I\{X_i \leq Y_i\}$ and $Y_i$ is independent of $X_i$.

***General partly interval censored data.*** The observed data are

$$\boldsymbol{O}_i = \begin{cases} X_i, & \text{if } 1 \leq i \leq n_1, \\ (\boldsymbol{Y}, \boldsymbol{\delta}_i), & \text{if } n_1 + 1 \leq i \leq n, \end{cases} \quad (1.7)$$

where for $N$ potential examination times $Y_1 < \cdots < Y_N$, letting $Y_0 = 0$ and $Y_{N+1} = \infty$, we have $\boldsymbol{Y} = (Y_1, \ldots, Y_N)$ and $\boldsymbol{\delta}_i = (\delta_i^{(1)}, \ldots, \delta_i^{(N+1)})$ with $\delta_i^{(j)} = 1$ if $Y_{j-1} < X_i \leq Y_j$ and 0 elsewhere. This means that for intervals $(0, Y_1], (Y_1, Y_2], \ldots, (Y_N, \infty)$, we know which one of them $X_i$ falls into. These two types of partial interval censoring were considered by Huang (1999), among others. As mentioned in Huang (1999), in practice, the general partly interval censored data were encountered in the Framingham Heart Disease Study (Odell, Anderson and D'Agostino (1992)) and in the study of the incidence of proteinuria in insulin-dependent diabetic patients (Enevoldsen *et al.* (1987)).

Obviously, one possible way to construct a likelihood-based confidence interval for $\theta_0$ with censored data is to use the likelihood function for a specific censoring model. This requires careful investigation for each type of censored sample. Specifically, the computation of the confidence region and the asymptotic results on the coverage of the confidence region need to be studied for each type of censored data. For works along these lines, still called the *empirical likelihood approach*, see Li, Hollander, McKeague and Yang (1996), Chen and Zhou (2003) and Banerjee and Wellner (2005) for right censored data, doubly censored data and interval censored Case 1 data, respectively. However, the methods in these works do not have direct extension to other types of censored data, such as interval censored Case 2 data (1.5) and partly interval censored data (1.6)–(1.7). Also, none of these works contains any coverage accuracy results or considers the use of smoothing. Note that the coverage accuracy of empirical likelihood ratio confidence intervals is $O(n^{-1})$ for smooth function models (DiCiccio, Hall and Romano (1991)). But, with censored data, we no longer have a smooth function model, thus it is quite difficult to study the coverage accuracy and it is particularly challenging to carry out the type of smoothing in Chen and Hall (1993) using the empirical likelihood approach for complicated types of censored data, such as doubly censored data, interval censored data and partly interval censored data.

Instead of undertaking the case-by-case studies for different types of censored data using the empirical likelihood approach, Ren (2001) constructs confidence intervals for the mean based on a new likelihood function, called a *weighted empirical likelihood function*, which is formulated in a unified form depending only on the probability mass of the *nonparametric maximum likelihood estimator* (NPMLE) $\hat{F}_n$ for $F_0$. For the mean, the $\sqrt{n}$-rate of convergence still holds for the aforementioned censored data (1.2)–(1.7) and Ren (2001) considered the first order expansion of the log-likelihood ratio without any coverage accuracy results. In this article, we construct smoothed *weighted empirical likelihood ratio confidence intervals* (WELRCI) for quantile $\theta_0$ in (1.1), where the $\sqrt{n}$-rate of convergence does not hold for interval censored data (1.4)–(1.5). Here, we derive the fourth order expansion of one type of smoothed weighted empirical log-likelihood ratio in a unified form for different types of censored data, including all of (1.2)–(1.7). With an analytically expressed leading term, this expansion leads to the *theoretical* coverage accuracy equation for the confidence intervals, which generally guarantees at least '*first order*' accuracy. This expansion also leads to the following results.

(a) When $\hat{F}_n$ has $\sqrt{n}$-rate of convergence, such as in the cases of right censored data, doubly censored data and partly interval censored data, the expansion shows that the log-likelihood ratio has an asymptotic scaled chi-squared distribution and the leading term of the expansion allows the $n$ out of $n$ bootstrap calibration for constructing confidence intervals in practice. Our theory shows that smoothed WELRCI is generally consistent; in particular, we show that for right censored data, the coverage accuracy of WELRCI is at least $O(n^{-1/2})$. Our simulation studies show that in comparison with empirical likelihood-based methods, the smoothing used in WELRCI generally gives a shorter confidence interval with comparable coverage accuracy.

(b) When $\hat{F}_n$ has $n^{1/3}$-rate of convergence, such as in the case of interval censored data, the expansion shows that with an adjusted rate of $n^{-1/3}$, the log-likelihood ratio has an asymptotic scaled $\mathbb{Z}^2$ distribution, where $n^{1/3}[\hat{F}_n(\theta_0) - F_0(\theta_0)] \xrightarrow{D} c_0\mathbb{Z}$; see (2.4) of Section 2 for details. It is interesting to notice that such a limiting distribution of the log-likelihood ratio is completely different from that obtained by the empirical likelihood approach for interval censored data (Banerjee and Wellner (2001, 2005)). In other words, the weighted empirical likelihood approach preserves the 'squared' structure of the limiting distribution of the log-likelihood ratio for all types of aforementioned censored data, while the empirical likelihood approach loses this structure for interval censored data. Moreover, for interval censored data, the aforementioned expansion of the weighted empirical log-likelihood ratio allows the use of the $m$ out $n$ bootstrap (Bickel, Götze and van Zwet (1997)) or subsampling (Politis, Romano and Wolf (1999)) calibration for constructing confidence intervals. With an adaptive choice of the bootstrap sample sizes which is similar to those by Götze and Račkauskas (2001) or Bickel and Sakov (2008), our simulation results for WELRCI using the $m$ out $n$ bootstrap compare favorably with those of Banerjee and Wellner (2005) for interval censored Case 1 data, though our method is computationally very time consuming for large $n$.

Regarding the main results of this article, an additional two points are worth mentioning.

(i) When there is no censoring, the bootstrap calibration of empirical likelihood has been studied by various authors (see the review in Owen (2001), Sections 3.3 and 3.17) and it is shown to have better coverage accuracy than standard chi-squared distribution calibration. However, for different types of censored data (1.2)–(1.7), it is not obvious how to generally implement the bootstrap calibration described in Owen (2001); for instance, it is not obvious how to correctly apply the $n$ out of $n$ bootstrap directly to the weighted empirical likelihood ratio, which itself is a solution of an optimization problem. Here, the expansion of the weighted empirical log-likelihood ratio provides a natural way to generally implement bootstrap calibration for censored data.

(ii) It is well known that the computation of likelihood-based confidence intervals is quite difficult in general. Here, the algorithm for computing WELRCI depends only on the NPMLE $\hat{F}_n$; that is, the routine itself is the same for different types of censored data. Thus, the WELRCI avoids complicated computation problems case-by-case.

The rest of this paper is organized as follows. Section 2 introduces the weighted empirical likelihood function with a brief review of the asymptotic properties of the NPMLE $\hat{F}_n$ for $F_0$. Section 3 constructs one type of smoothed WELRCI for quantile $\theta_0$ and gives related asymptotic results with proofs deferred to Section 6. Section 4 discusses the implementation of WELRCI in practice. Section 5 presents some simulation results and comparisons between the proposed procedure and alternative methods, and includes some concluding remarks.

## 2. Weighted empirical likelihood

In Owen (1988), the empirical likelihood function is given by

$$L(F) = \prod_{i=1}^{n}[F(X_i) - F(X_{i-})], \qquad (2.1)$$

where $F$ is any d.f., and the empirical likelihood ratio function is given by $R(F) = L(F)/L(F_n)$ because the empirical d.f. $F_n$ of sample $X_1, \ldots, X_n$ is the *nonparametric maximum likelihood estimator* (NPMLE) for $F_0$; that is, $F_n$ maximizes $L(F)$ over all distribution functions $F$. The weighted empirical likelihood function in Ren (2001) is given as follows.

For each type of aforementioned censored data (1.2)–(1.7), the NPMLE $\hat{F}_n$ for $F_0$ based on observed data can be expressed as

$$\hat{F}_n(x) = \sum_{i=1}^{m} \hat{p}_i I\{W_i \le x\}, \qquad (2.2)$$

where $W_1 < W_2 < \cdots < W_m$ with $\hat{p}_j > 0$, $1 \le j \le m$. Specifically, in the case of right censored data, $W_i$'s are those non-censored observations in (1.2) and $m$ is the total number of non-censored observations among $V_1 < \cdots < V_n$; see Kaplan and Meier (1958) or Shorack and Wellner (1986), page 293. In the case of doubly censored data, $W_i$'s include all those non-censored observations among $V_1 < \cdots < V_n$ in (1.3), but certain $W_j$ could be a $V_k$ with, say, $\delta_k = 3$; see Mykland and Ren (1996). For doubly censored data, the NPMLE $\hat{F}_n$ is implicitly, but uniquely, determined by a particular solution of an integral equation; in turn, $W_i$'s and $\hat{p}_i$'s are uniquely determined and can be obtained through computing $\hat{F}_n$ (Mykland and Ren (1996)). In the case of interval censored Case 2 data (1.5), the NPMLE $\hat{F}_n$ is also implicitly, but uniquely, determined by a particular solution of an integral equation; in turn, $W_i$'s and $\hat{p}_i$'s are uniquely determined and can be obtained through computing $\hat{F}_n$; see Groeneboom and Wellner (1992) pages 43–46. For interval censored Case 2 data, $\{W_1, \ldots, W_m\}$ is a subset of $\{Y_1, \ldots, Y_n, Z_1, \ldots, Z_n\}$ in (1.5). Similarly, for interval censored Case 1 data (1.4) and partly interval censored data (1.6)–(1.7), $W_i$'s and $\hat{p}_i$'s in (2.2) are also uniquely determined by computing $\hat{F}_n$; see Groeneboom and Wellner (1992), pages 35–41 and Huang (1999).

The *weighted empirical likelihood function* (Ren (2001)) is given by

$$\hat{L}(F) = \prod_{i=1}^{m} [F(W_i) - F(W_{i^-})]^{n\hat{p}_i} \qquad (2.3)$$

and it is shown (Ren (2008)) that $\hat{L}(F)$ may be viewed as the asymptotic version of the empirical likelihood function $L(F)$ for censored data. It is easy to show that $\hat{L}(F)$ is maximized at $\hat{F}_n$. Hence, the *weighted empirical likelihood ratio* is given by $\hat{R}(F) = \hat{L}(F)/\hat{L}(\hat{F}_n)$. One may notice that when there is no censoring, the weighted empirical likelihood function (2.3) coincides with Owen's empirical likelihood function (2.1); see Ren (2001) for details.

**Remark 1.** *Asymptotic results for the NPMLE $\hat{F}_n$.* It is known that $\|\hat{F}_n - F_0\| \xrightarrow{a.s.} 0$ as $n \to \infty$ for right censored data (Stute and Wang (1993)), doubly censored data (Gu and Zhang (1993)), interval censored data (Groeneboom and Wellner (1992)), and partly interval censored data (Huang (1999)), respectively. It is also known that for right censored data, doubly censored data and partly interval censored data, $\sqrt{n}(\hat{F}_n - F_0)$ weakly converges to a centered Gaussian process under certain conditions (Gill (1983), Gu and Zhang (1993), Huang (1999)). However, for interval censored Case 1 data (1.4), we have that for a fixed point $t_0$,

$$n^{1/3}[\hat{F}_n(t_0) - F_0(t_0)] \xrightarrow{D} c_0 \mathbb{Z} \qquad \text{as } n \to \infty, \qquad (2.4)$$

where $c_0$ is a constant and $\mathbb{Z} = \arg\min(W(t) + t^2)$ with $W$ being the two-sided Brownian motion starting from 0 (Groeneboom and Wellner (1992)). For interval censored Case 2 data (1.5), (2.4) also holds under certain conditions (Wellner (1995)). Note that (2.4) accents for why a $\sqrt{n}$-rate of convergence does not hold for quantile estimators with interval censored data.

## 3. Weighted empirical likelihood ratio confidence intervals

In this section, we show that the set $S_n = \{\tilde{F}^{-1}(q) | \hat{R}(F) \geq c_n, F \ll \hat{F}_n\}$ may be used as confidence interval for the quantile $\theta_0$ given in (1.1), where $0 < c_n < 1$ is a constant, $\tilde{F}$ is a smoothed version of $F$, as given in equation (3.2) below, and '$F \ll \hat{F}_n$' means that $F$ is absolutely continuous with respect to $\hat{F}_n$.

First, note that the NPMLE $\hat{F}_n$ for censored data (1.2)–(1.7) is not always a proper d.f. (Mykland and Ren (1996)), but, in this work, we always consider the adjusted version of the NPMLE, still denoted $\hat{F}_n$. Precisely, for the rest of this paper, $\hat{F}_n$ in (2.2) denotes the proper d.f. obtained by setting 1 as the value of the NPMLE at the largest observation of the data set, that is, setting $\hat{F}_n = 1$ at $V_{(n)}, Y_{(n)}, \max\{Y_{(n)}, Z_{(n)}\}$ or $\max\{X_i\text{'s}, Y_j\text{'s}\}$, which implies that $\sum_{i=1}^{m} \hat{p}_i = 1$ in (2.2). This kind of adjustment of the NPMLE is a generally adopted convention for censored data (Efron (1967); Miller (1976)). Although this $\hat{F}_n$ no longer necessarily maximizes the underlying likelihood function, the usual asymptotic properties of the NPMLE needed for this work still hold for this $\hat{F}_n$; see the later discussion in Remark 2.

To study the confidence set $S_n$, we let $p_i = F(W_i) - F(W_{i^-}), 1 \leq i \leq m$, and let

$$r(\theta) = \sup\left\{\prod_{i=1}^{m}(p_i/\hat{p}_i)^{n\hat{p}_i} \,\Big|\, \tilde{F}_{\boldsymbol{p}}^{-1}(q) = \theta, F_{\boldsymbol{p}} \in \mathfrak{F}_n\right\}, \tag{3.1}$$

where $\mathfrak{F}_n \equiv \{F \mid F(x) = \sum_{i=1}^{m} p_i I\{W_i \leq x\}, p_i \geq 0, \sum_{i=1}^{m} p_i = 1\}$, and for $W_0 = 0, \boldsymbol{W} = (W_1, \ldots, W_m)$ and $\boldsymbol{p} = (p_1, \ldots, p_m)$, $\tilde{F}_{\boldsymbol{p}}$ is a smoothed version of $F_{\boldsymbol{p}}$ by connecting adjacent jump points through straight lines for $0 < x \leq W_m$:

$$\tilde{F}_{\boldsymbol{p}}(x) = \sum_{i=1}^{m} I\{W_{i-1} < x \leq W_i\}\left(\sum_{j=1}^{i-1} p_j + \frac{p_i(x - W_{i-1})}{W_i - W_{i-1}}\right) = \sum_{i=1}^{m} p_i H_i(\boldsymbol{W}, x) \tag{3.2}$$

with $H_i(\boldsymbol{W}, x) = \frac{x - W_{i-1}}{W_i - W_{i-1}} I\{W_{i-1} < x \leq W_i\} + I\{W_i < x\}$. The proof of the last equation in (3.2) is based on straightforward algebra, which is omitted for brevity. In the Appendix, we show that $S_n$ is an interval satisfying $S_n = [X_L, X_U]$ and

$$X_L \leq \theta_0 \leq X_U \quad \text{if and only if} \quad r(\theta_0) \geq c_n, \tag{3.3}$$

where

$$\begin{aligned} X_L &= \inf\left\{\tilde{F}_{\boldsymbol{p}}^{-1}(q) \,\Big|\, F_{\boldsymbol{p}} \in \mathfrak{F}_n, \prod_{i=1}^{m}(p_i/\hat{p}_i)^{n\hat{p}_i} \geq c_n\right\}, \\ X_U &= \sup\left\{\tilde{F}_{\boldsymbol{p}}^{-1}(q) \,\Big|\, F_{\boldsymbol{p}} \in \mathfrak{F}_n, \prod_{i=1}^{m}(p_i/\hat{p}_i)^{n\hat{p}_i} \geq c_n\right\}. \end{aligned} \tag{3.4}$$

We call $[X_L, X_U]$ the smoothed *weighted empirical likelihood ratio confidence interval* (WELRCI) for the quantile $\theta_0$. Note that (3.3) does not hold if, in (3.1), we use $F_{\boldsymbol{p}}^{-1}(q) = \theta$ in place of $\tilde{F}_{\boldsymbol{p}}^{-1}(q) = \theta$ because $F_{\boldsymbol{p}}^{-1}(q) = \theta$ is not equivalent to $F_{\boldsymbol{p}}(\theta) = q$. Also, note that other types of smoothing, such as the kernel density estimator method (Chen and Hall (1993), Ren (2006)), may be considered in (3.1). The smoothing issue will be discussed later in Remark 4.

Since (3.3) implies

$$P\{X_L \leq \theta_0 \leq X_U\} = P\{-2\log r(\theta_0) \leq -2\log c_n\}, \tag{3.5}$$

the asymptotic behavior of $[X_L, X_U]$ is studied via the weighted empirical log-likelihood ratio $\log r(\theta_0)$ in Theorem 1 with proofs given in Section 6. In Theorem 1, we let

$$\hat{\theta} = \tilde{F}_n^{-1}(q) \quad \text{and} \quad \hat{\eta} = \tilde{F}_n(\theta_0), \tag{3.6}$$

where $\tilde{F}_n$ denotes the smoothed version of $\hat{F}_n$ according to (3.2) and we let

$$\hat{\mu}_k = \hat{\mu}_k(\theta_0) \quad \text{and} \quad \hat{\mu}_k(\theta) \equiv \sum_{i=1}^{m} \hat{p}_i [H_i(\boldsymbol{W}, \theta) - q]^k, \qquad k = 1, 2, \ldots. \tag{3.7}$$

**Theorem 1.** *Assume that for a sequence $C_n \to \infty$, we have that as $n \to \infty$,*

$$C_n(\hat{\eta} - q) = O_p(1), \tag{AS1}$$

$$\hat{\eta} \xrightarrow{a.s.} q, \tag{AS2}$$

$$\hat{\mu}_2 \xrightarrow{a.s.} q(1-q). \tag{AS3}$$

*Then:*

(i) *with probability 1, we have that for fixed $k = 0, 1, 2, 3, 4$,*

$$-2\log r(\theta_0) = B_n^{(k)} + n(\hat{\eta} - q)^{k+3} r_{n,k}, \qquad |r_{n,k}| \leq M_{r,k}, \tag{3.8}$$

*all but finitely often, where $1 \leq M_{r,k} < \infty$ is a constant and*

$$B_n^{(k)} = \frac{n(\hat{\eta} - q)^2}{\hat{\mu}_2} \left(1 + \sum_{j=1}^{k} \hat{a}_j (\hat{\eta} - q)^j\right) \tag{3.9}$$

*with $B_n^{(0)} = n(\hat{\eta} - q)^2 / \hat{\mu}_2$ and*

$$\begin{aligned}
\hat{a}_1 &= (2\hat{\mu}_3)/(3\hat{\mu}_2^2), \qquad \hat{a}_2 = \hat{\mu}_2^{-4}(\hat{\mu}_3^2 - \tfrac{1}{2}\hat{\mu}_2 \hat{\mu}_4), \\
\hat{a}_3 &= 2\hat{\mu}_2^{-6}(\hat{\mu}_3^3 + \tfrac{1}{5}\hat{\mu}_2^2 \hat{\mu}_5 - \hat{\mu}_2 \hat{\mu}_3 \hat{\mu}_4), \\
\hat{a}_4 &= \hat{\mu}_2^{-8}(\tfrac{14}{3}\hat{\mu}_3^4 - \tfrac{1}{3}\hat{\mu}_2^3 \hat{\mu}_6 + \hat{\mu}_2^2 \hat{\mu}_4^2 + 2\hat{\mu}_2^2 \hat{\mu}_3 \hat{\mu}_5 - 7\hat{\mu}_2 \hat{\mu}_3^2 \hat{\mu}_4);
\end{aligned} \tag{3.10}$$

(ii) *assuming that $c_n$ is chosen such that $\tilde{c}_n = n^{-1}(C_n)^2(-2\log c_n) = O(1)$ and assuming that*

$$[C_n(\hat{\eta} - q)]^2/\hat{\mu}_2 \text{ has a limiting distribution } G_0, \qquad (AS4)$$

*where $G_0$ is continuous with bounded derivative on some finite interval $[a,b]$, which satisfies $[\tilde{c}_n - \delta, \tilde{c}_n + \delta] \subset [a,b]$ for some $a > 0, \delta > 0$ and sufficiently large $n$, then,*

$$P\{X_L \leq \theta_0 \leq X_U\} = P\{A_n^{(k)} \leq \tilde{c}_n\} + O(\|F_{n,k} - G_0\|_{[a,b]}) + O((C_n)^{-(k+1)}), \qquad (3.11)$$

*where $F_{n,k}$ is the d.f. of $A_n^{(k)} = n^{-1}(C_n)^2 B_n^{(k)}$ and $\|\cdot\|_{[a,b]}$ is the uniform norm on $[a,b]$.*

In practice, we let $\rho_{n,\alpha}^{(k)}$ be the $(1-\alpha)100$th percentile of $A_n^{(k)}$ in (3.11) for $0 < \alpha < 1$ and then $[X_L^{(k)}, X_U^{(k)}]$ computed by (3.4) with constant $c_n$ set by

$$-2\log c_n = n(C_n)^{-2}\rho_{n,\alpha}^{(k)} \quad \Leftrightarrow \quad \tilde{c}_n = \rho_{n,\alpha}^{(k)} \qquad (3.12)$$

is called the *kth order WELRCI* (*k*-WELRCI) for $\theta_0$. Thus, from (3.11)–(3.12), we have the *theoretical coverage accuracy equation*

$$P\{X_L^{(k)} \leq \theta_0 \leq X_U^{(k)}\} = (1-\alpha) + O(\|F_{n,k} - G_0\|_{[a,b]}) + O((C_n)^{-(k+1)}), \qquad (3.13)$$

where the convergence rate of $\|F_{n,k} - G_0\|_{[a,b]}$ is referred to as the '*first order*'. In Remark 5 and Section 4, the coverage accuracy issue of *k*-WELRCI and the estimation of $\rho_{n,\alpha}^{(k)}$ in practice will be discussed, respectively.

As mentioned in Section 1, for censored data, we no longer have a smooth function model and at present, the coverage accuracy of likelihood-based confidence intervals is unknown. With the proof deferred to Section 6, the next corollary shows that (3.8) and (3.11) can help us study the coverage accuracy of *k*-WELRCI.

**Corollary 1.** *Under the assumptions of Theorem 1 and the assumption in (2.1) of Chen and Lo (1996), the coverage accuracy of WELRCI for the quantile $\theta_0$ with right censored data is at least $O(n^{-1/2})$, that is, with $k = 0$ in (3.11) and $c_n$ set by (3.12), we have*

$$P\{X_L^{(0)} \leq \theta_0 \leq X_U^{(0)}\} = (1-\alpha) + O(n^{-1/2}). \qquad (3.14)$$

**Remark 2.** *Assumptions of Theorem* 1. Since $\tilde{F}_n$ is an increasing function on the support of $\hat{F}_n$ and since

$$\|\tilde{F}_n - \hat{F}_n\| \leq \sup_x |\hat{F}_n(x) - \hat{F}_{n(x^-)}| = O(\|\hat{F}_n - F_0\|), \qquad (3.15)$$

we have $q = F_0(\theta_0) = \tilde{F}_n(\hat{\theta})$ and that the asymptotic properties of $(\hat{\eta} - q) = [\tilde{F}_n(\theta_0) - F_0(\theta_0)]$ are determined by those of $[\hat{F}_n(\theta_0) - F_0(\theta_0)]$. From Remark 1, we know that under suitable conditions, (AS2) holds for all those censored data (1.2)–(1.7). Also, it is easy to show that $\|\hat{F}_n - F_0\| \xrightarrow{a.s.} 0$ implies (AS3). For (AS1), note that whenever it is known, $C_n$ is meant to be the convergence rate of $\hat{\eta}$. Thus, for right censored data, doubly censored data and partly interval censored data, from Remark 1, we know that $\sqrt{n}[\hat{F}_n(\theta_0) - F_0(\theta_0)]$ has an asymptotic normal distribution; in turn, under suitable conditions, (AS1) and (AS4) hold with $C_n = \sqrt{n}$. Moreover, for right censored data and doubly censored data, we have $|\hat{F}_n(\theta_0) - \tilde{F}_n(\theta_0)| = O_{a.s.}(n^{-1})$ (Ren (1997)), which implies that $\sqrt{n}[\hat{F}_n(\theta_0) - F_0(\theta_0)]$ and $\sqrt{n}[\tilde{F}_n(\theta_0) - F_0(\theta_0)]$ have the same limiting distribution; in turn, $G_0$ is the d.f. of $\rho_0 \chi_1^2$, where $\rho_0$ is some constant and $\chi_1^2$ is a chi-squared random variable with one degree of freedom. For interval censored Case 1 or Case 2 data, we have (2.4) under certain conditions, and from Groeneboom and Wellner (1992), we could expect $n^{1/3}[\hat{F}_n(W_j) - \hat{F}_n(W_j-)] = o_p(1)$ for $W_{j-1} \leq \theta_0 \leq W_j$, thus (AS1) and (AS4) hold with $C_n = n^{1/3}$ when (2.4) holds, where $G_0$ is the d.f. of $\gamma_0 \mathbb{Z}^2$ for some constant $\gamma_0$. On the other hand, if smoothing is not used in (3.1), that is, replacing $\tilde{F}_{\boldsymbol{p}}^{-1}(q) = \theta$ by $F_{\boldsymbol{p}}(\theta) = q$ in (3.1), a modified proof of Theorem 1 shows that (3.8)–(3.11) hold with $\hat{\eta} = \hat{F}_n(\theta_0)$, for which (AS1) and (AS4) hold with $C_n = n^{1/3}$ for interval censored data whenever (2.4) holds. In our simulation studies presented in Section 5, we denote this non-smoothed WELRCI as WELRCI0. Finally, note that (3.8)–(3.9), (AS3)–(AS4), the continuity of $G_0$ and Pólya's theorem imply that $\|F_{n,k} - G_0\|_{[a,b]} \to 0$ as $n \to \infty$. Thus, if $c_n$ is set by (3.12), we have $\tilde{c}_n = O(1)$ in (3.11).

***Remark 3.*** *Limiting distribution of log-likelihood ratio.* When $\hat{F}_n$ has $\sqrt{n}$-rate of convergence, such as in the cases of right censored data, doubly censored data and partly interval censored data, from (3.8)–(3.9), $|\hat{a}_k| \leq 1$, (AS3) and Remark 2, we know that the weighted empirical log-likelihood ratio $-2 \log r(\theta_0) \xrightarrow{D} \rho_0 \chi_1^2$, where $\rho_0 = 1$ when there is no censoring. Similarly, when $\hat{F}_n$ has $n^{1/3}$-rate of convergence, such as in the case of interval censored data, we know that (3.8)–(3.9), (2.4) and Remark 2 give $n^{-1/3}[-2 \log r(\theta_0)] \xrightarrow{D} \gamma_0 \mathbb{Z}^2$, while the empirical log-likelihood ratio of Banerjee and Wellner (2001), page 1701; (2005), page 411 converges in distribution to $\mathbb{D}$, which is not proportional to $\mathbb{Z}^2$.

***Remark 4.*** *Smoothing and theoretical coverage accuracy of WELRCI.* From (6.21) in the proof of Theorem 1(ii), it is easy to see that if $F_{n,k}$ satisfies the Lipschitz condition $|F_{n,k}(x) - F_{n,k}(y)| \leq M_F |x-y|$ in the neighborhood of $\tilde{c}_n$ for some constant $M_F$ and sufficiently large $n$, then the term $O(\|F_{n,k} - G_0\|_{[a,b]})$ in equation (3.13) disappears, that is, (3.13) becomes the following *best possible theoretical coverage accuracy equation:*

$$P\{X_L^{(k)} \leq \theta_0 \leq X_U^{(k)}\} = (1-\alpha) + O((C_n)^{-(k+1)}). \quad (3.16)$$

Note that the use of $\tilde{F}_{\boldsymbol{p}}$ in (3.1) leads to $\hat{\eta} = \tilde{F}_n(\theta_0)$ in (3.6) and it can be shown that without smoothing, that is, using $F_{\boldsymbol{p}}(\theta) = q$ in place of $\tilde{F}_{\boldsymbol{p}}^{-1}(q) = \theta$ in (3.1), Theorem 1 has $\hat{\eta} = \hat{F}_n(\theta_0)$, for which the Lipschitz condition does not hold. On the other hand, in Ren (2006), it is shown that if $\tilde{F}_n$ is based on the kernel density method, then the Lipschitz condition for $F_{n,k}$ holds asymptotically when there is no censoring. The implication of this is that it is possible to have (3.16) for censored data with an appropriate smoothing of $\hat{F}_n$. Due to this understanding on the effects of smoothing $\hat{F}_n$, we consider a very simple smoothed version $\tilde{F}_p$ of $F_p$ in (3.1) and the simulation results compare favorably with alternative methods. However, in general, the verification of the Lipschitz condition can be quite involved for censored data; this will be studied in a separate paper. Thus, for now, we may only view (3.16) as the *theoretically* best possible for the $k$-WELRCI, which will be used in the next section to select $k$ in practice. It should be noted that (3.16) suggests that the use of the fourth order expansion in (3.8) is sufficient because $k = 4$ and $C_n = \sqrt{n}$ give $O(n^{-5/2})$ for the last term in (3.16), while the coverage accuracy with Bartlett-correction is only $O(n^{-2})$ for smooth function models (DiCiccio, Hall and Romano (1991)).

**Remark 5.** *Coverage accuracy.* The coverage accuracy equation (3.13) is only *theoretical*, because the actual coverage accuracy in practice includes the estimation error for $\rho_{n,\alpha}^{(k)}$. This means that the actual coverage accuracy can be established via the rate of $\|F_{n,k} - G_0\|_{[a,b]}$ and the estimation error of $\rho_{n,\alpha}^{(k)}$. But, if we have (3.16) via an appropriate smoothing, the coverage accuracy is determined only through the estimation error of $\rho_{n,\alpha}^{(k)}$. In the case of right censored data, our (3.14) in Corollary 1, along with the bootstrap estimation error results established in Chen and Lo (1996), ensures the actual coverage accuracy of our WELRCI to be at least $O(n^{-1/2})$. For other types of censored data (1.3)–(1.7), our equations (3.12)–(3.13) and (3.16) indicate the direction of further studies on the actual coverage accuracy and provide guidance on the implementation of our smoothed WELRCI in practice, which is discussed in the next section.

## 4. Implementation

### Estimation of $\rho_{n,\alpha}^{(k)}$

Since the computation of $[X_L, X_U]$ in (3.4) only depends on the constant $c_n$, (3.12) implies that $k$-WELRCI $[X_L^{(k)}, X_U^{(k)}]$ for the quantile $\theta_0$ can be computed if $\rho_{n,\alpha}^{(k)}$ can be estimated consistently for an appropriate $k$.

For right censored data, doubly censored data and partly interval censored data, we have $C_n = \sqrt{n}$ and, in turn, $A_n^{(k)} = B_n^{(k)}$ in (3.11) can be expressed as

$$A_n^{(k)} = \frac{[C_n(\hat{\eta} - q)]^2}{\hat{\mu}_2}\left(1 + \sum_{j=1}^{k} \hat{b}_j [C_n(\hat{\eta} - q)]^j\right) \equiv \tau_n(C_n(\hat{\eta} - q)), \tag{4.1}$$

where $\hat{b}_j = \hat{a}_j/(C_n)^j$. Thus, $\rho_{n,\alpha}^{(k)}$ may be estimated by the percentiles of

$$A_n^{(k)*} = \frac{[C_n(\hat{\eta}^* - q)]^2}{\hat{\mu}_2(\hat{\theta})}\left(1 + \sum_{j=1}^{k} \hat{b}_j(\hat{\theta})[C_n(\hat{\eta}^* - q)]^j\right) \equiv \hat{\tau}_n(C_n(\hat{\eta}^* - q)) \tag{4.2}$$

for $\hat{b}_j(\hat{\theta}) = \hat{a}_j(\hat{\theta})/(C_n)^j$ and $\hat{\eta}^* = \tilde{F}_n^*(\hat{\theta})$, where for $\hat{\mu}_k(\hat{\theta})$ given by (3.6)–(3.7), $\hat{a}_j(\hat{\theta})$ are calculated by (3.10) with $\hat{\mu}_k$'s replaced by $\hat{\mu}_k(\hat{\theta})$'s and $\tilde{F}_n^*$ is calculated based on the $n$ out of $n$ bootstrap method (Efron (1979)). That is, if $\rho_{n,\alpha}^{(k)*}$ denotes the $(1-\alpha)100$th percentile of $A_n^{(k)*}$, we use it to estimate $\rho_{n,\alpha}^{(k)}$ in practice. Note that the asymptotic properties of $(\hat{\eta} - q) = [\tilde{F}_n(\theta_0) - F_0(\theta_0)]$ are the same as those of $[\hat{F}_n(\theta_0) - F_0(\theta_0)]$ due to the smoothing method used in (3.2). Thus, from Bickel and Ren (1996) and Huang (1999), the bootstrap consistency holds here because it is easy to show that $\|\hat{F}_n - F_0\| \xrightarrow{a.s.} 0$ implies $\hat{\theta} \xrightarrow{a.s.} \theta_0$ and $|\hat{\mu}_k(\hat{\theta}) - \hat{\mu}_k| \xrightarrow{a.s.} 0$; in turn, $|\hat{\alpha}_j(\hat{\theta}) - \hat{\alpha}_j| \xrightarrow{a.s.} 0$ as $n \to \infty$.

For interval censored data, if (2.4) holds, we have $C_n = n^{1/3}$ in (4.1) and the distribution of $A_n^{(k)} = n^{-1/3} B_n^{(k)}$ may be estimated by that of $\hat{\tau}_n(n_b^{1/3}(\hat{\eta}_{n_b}^* - q))$, where $\hat{\eta}_{n_b}^* = \tilde{F}_{n_b}^*(\hat{\theta})$ is calculated based on the subsampling method (Politis, Romano and Wolf (1999), Theorem 2.2.1) or the $m$ out of $n$ bootstrap method (Bickel, Götze and van Zwet (1997)), using $n_b$ as the resampling size. In this case, $\rho_{n,\alpha}^{(k)}$ is estimated by the $(1-\alpha)100$th percentile of $\hat{\tau}_n(n_b^{1/3}(\hat{\eta}_{n_b}^* - q))$. With an adaptive choice of the bootstrap sample size $n_b$, the $m$ out of $n$ bootstrap performs very well in our simulation studies, as shown in Section 5.

### Selection of $k$

Note that Corollary 1, along with Remark 5, shows that the coverage accuracy of WEL-RCI is at least as good as that achieved by normal-based methods for right censored data. But, this may not be the exact coverage accuracy. In fact, when there is no censoring, the term $O(\|F_{n,k} - G_0\|_{[a,b]})$ (which usually has rate $n^{-1/2}$) in (3.13) does not exist because the coverage accuracy is $O(n^{-1})$ for smooth function models (DiCiccio, Hall and Romano (1991)). Thus, based on Remark 4, we know that, at best, we may expect the theoretical coverage accuracy equation (3.16) to hold for $k$-WELRCI, which may be used to select $k$ in practice. Since the estimation accuracy for $\rho_{n,\alpha}^{(k)}$ via the usual bootstrap method cannot be better than $O_p(n^{-1})$ (Hall (1992)), the criterion we recommend for selecting $k$ in practice is to choose the smallest $k$ such that

$$(C_n)^{-(k+1)} < n^{-1}, \tag{4.3}$$

which implies that $n(C_n)^{-(k+3)} = o(1)$ in (3.8). Thus, for right censored data, doubly censored data and partly interval censored data, we have $C_n = \sqrt{n}$, for which (4.3) implies $k = 2$; for interval censored Case 1 data (1.4), we have $C_n = n^{1/3}$, for which (4.3) implies $k = 3$.

### Computation

A routine in FORTRAN for computing $[X_L^{(k)}, X_U^{(k)}]$ based on (3.4) is available from the author, but for the brevity of this article, the theorems on the convergence of this routine are omitted.

**Table 1.** 90% confidence intervals for $\theta_0 = f_0^{-1}(q)$ with exponential right censored data $X \sim \text{Exp}(1)$, $Y \sim \text{Exp}(3)$; percentage of $\delta$: $\delta = 1$: 75.0%; $\delta = 0$: 25.0%

| | | $n = 50$ | | $n = 100$ | | $n = 200$ | |
|---|---|---|---|---|---|---|---|
| $q$ & $\theta_0$ | Method | Coverage % | Average length (s.d.) | Coverage % | Average length (s.d.) | Coverage % | Average length (s.d.) |
| $q = 0.250$ | 1-WELRCI | 90.1 | 0.272 (0.090) | 89.9 | 0.191 (0.051) | 89.8 | 0.136 (0.032) |
| $\theta_0 = 0.288$ | 2-WELRCI | 90.5 | 0.273 (0.090) | 89.9 | 0.192 (0.051) | 89.8 | 0.136 (0.032) |
| | LHMYCI | 90.1 | 0.285 (0.095) | 90.3 | 0.196 (0.052) | 90.2 | 0.138 (0.032) |
| | SQBPCI | 90.4 | 0.267 (0.086) | 90.6 | 0.190 (0.049) | 90.0 | 0.136 (0.032) |
| | QBPCI | 89.4 | 0.286 (0.098) | 90.0 | 0.196 (0.053) | 89.5 | 0.138 (0.033) |
| $q = 0.500$ | 1-WELRCI | 89.8 | 0.480 (0.151) | 91.0 | 0.344 (0.091) | 90.0 | 0.242 (0.057) |
| $\theta_0 = 0.693$ | 2-WELRCI | 89.4 | 0.473 (0.150) | 90.7 | 0.342 (0.091) | 90.0 | 0.241 (0.057) |
| | LHMYCI | 89.9 | 0.513 (0.167) | 91.8 | 0.357 (0.094) | 90.0 | 0.246 (0.057) |
| | SQBPCI | 88.5 | 0.467 (0.147) | 90.6 | 0.341 (0.089) | 89.4 | 0.241 (0.056) |
| | QBPCI | 89.2 | 0.506 (0.169) | 91.6 | 0.354 (0.095) | 89.7 | 0.245 (0.058) |
| $q = 0.750$ | 1-WELRCI | 88.6 | 0.860 (0.329) | 88.4 | 0.637 (0.211) | 90.0 | 0.460 (0.128) |
| $\theta_0 = 1.386$ | 2-WELRCI | 88.8 | 0.866 (0.331) | 88.7 | 0.639 (0.211) | 90.0 | 0.461 (0.128) |
| | LHMYCI | 90.5 | 1.057 (0.512) | 90.0 | 0.706 (0.251) | 91.3 | 0.484 (0.136) |
| | SQBPCI | 87.2 | 0.832 (0.297) | 88.0 | 0.624 (0.201) | 89.4 | 0.457 (0.123) |
| | QBPCI | 90.7 | 1.041 (0.499) | 89.1 | 0.684 (0.243) | 91.2 | 0.476 (0.136) |

**Table 2.** 90% confidence intervals for $\theta_0 = F_0^{-1}(q)$ with chi-squared right censored data $X \sim \chi^2(1)$, $Y \sim \text{Exp}(3)$; percentage of $\delta$: $\delta = 1$: 77.5%; $\delta = 0$: 22.5%

| | | $n = 50$ | | $n = 100$ | | $n = 200$ | |
|---|---|---|---|---|---|---|---|
| $q$ & $\theta_0$ | Method | Coverage % | Average length (s.d.) | Coverage % | Average length (s.d.) | Coverage % | Average length (s.d.) |
| $q = 0.250$ | 1-WELRCI | 89.5 | 0.172 (0.075) | 89.5 | 0.123 (0.044) | 90.0 | 0.087 (0.026) |
| $\theta_0 = 0.102$ | 2-WELRCI | 89.8 | 0.173 (0.076) | 89.7 | 0.124 (0.044) | 90.0 | 0.087 (0.026) |
| | LHMYCI | 90.4 | 0.187 (0.083) | 88.9 | 0.128 (0.045) | 90.1 | 0.089 (0.026) |
| | SQBPCI | 89.7 | 0.169 (0.074) | 89.4 | 0.123 (0.043) | 90.1 | 0.087 (0.026) |
| | QBPCI | 88.8 | 0.191 (0.087) | 89.0 | 0.128 (0.047) | 90.0 | 0.090 (0.027) |
| $q = 0.500$ | 1-WELRCI | 89.8 | 0.500 (0.200) | 88.7 | 0.359 (0.110) | 90.1 | 0.253 (0.063) |
| $\theta_0 = 0.455$ | 2-WELRCI | 89.2 | 0.492 (0.197) | 88.4 | 0.356 (0.109) | 89.9 | 0.252 (0.063) |
| | LHMYCI | 89.9 | 0.546 (0.231) | 88.9 | 0.374 (0.117) | 89.4 | 0.259 (0.064) |
| | SQBPCI | 88.4 | 0.484 (0.193) | 89.0 | 0.354 (0.108) | 90.1 | 0.252 (0.062) |
| | QBPCI | 89.1 | 0.542 (0.230) | 88.8 | 0.371 (0.118) | 89.8 | 0.258 (0.065) |
| $q = 0.750$ | 1-WELRCI | 85.9 | 1.131 (0.465) | 87.1 | 0.847 (0.308) | 88.4 | 0.630 (0.184) |
| $\theta_0 = 1.323$ | 2-WELRCI | 86.5 | 1.138 (0.467) | 87.1 | 0.850 (0.308) | 88.4 | 0.631 (0.184) |
| | LHMYCI | 85.5 | 1.438 (0.871) | 90.0 | 0.963 (0.403) | 89.6 | 0.669 (0.202) |
| | SQBPCI | 84.4 | 1.083 (0.428) | 86.1 | 0.831 (0.287) | 88.2 | 0.625 (0.178) |
| | QBPCI | 86.8 | 1.470 (0.830) | 89.0 | 0.943 (0.402) | 89.6 | 0.660 (0.201) |

## 5. Simulation studies and concluding remarks

Some simulation results for Theorem 1 are presented in this section. Here, letting $\text{Exp}(\mu)$ and $\chi^2(r)$ represent the exponential d.f. with mean $\mu$ and the chi-squared d.f. with degrees of freedom $r$, respectively, the order $k$ for smoothed $k$-WELRCI is selected based on (4.3) and $[X_L^{(k)}, X_U^{(k)}]$ is computed by the algorithm mentioned in Section 4. In all of our simulation studies here, the EM algorithm is used to compute the NPMLE $\hat{F}_n$ for doubly censored data and interval censored data, and the stopping rule used is when the uniform distance between two consecutive iterations is less than 0.001.

### Right censored data and doubly censored data

Since the empirical likelihood-based confidence intervals (C.I.) for quantiles were considered by Li, Hollander, McKeague and Yang (1996) for right censored data (abbreviated as LHMYCI), we make comparisons between WELRCI and LHMYCI. Moreover, other procedures, such as bootstrap percentile confidence intervals (Efron and Tibshirani (1993)) for sample quantiles $\hat{F}_n^{-1}(q)$ (abbreviated as QBPCI) and for smoothed quantiles $\hat{\theta} = \tilde{F}_n^{-1}(q)$ in (3.6) (abbreviated as SQBPCI), are also considered in our studies. In Table 1, 1000 right censored samples (1.2) of size $n = 50$ are taken from exponential distributions and for each sample, 90% $k$-WELRCI, LHMYCI, SQBPCI and QBPCI for $\theta_0 = F_0^{-1}(q)$ with different $q$ are computed, where 400 bootstrap samples of size $n = 50$ are used for SQBPCI, QBPCI, and for estimating $\rho_{n,\alpha}^{(k)}$ in (3.12) to construct smoothed $k$-WELRCI. The simulation coverage is included in Table 1, and the simulation standard deviation (s.d.) of the length of C.I. is given in the parentheses next to the average length of C.I. Table 1 also includes the results of the same studies with $n = 100$ and $n = 200$, respectively. The simulation studies in Table 1 are repeated in Table 2 with $F_0 = \chi^2(1)$ and repeated in Tables 3 and 4 with exponential and chi-squared doubly censored samples (1.3), respectively.

From Tables 1–2, we see that for right censored data, WELRCI and LHMYCI have comparable coverage accuracy, but LHMYCI, while not always having better coverage, are noticeably wider than WELRCI for moderate sample size, say, $n = 50$. For right censored data and doubly censored data, WELRCI have coverage accuracy similar to QBPCI, but are noticeably shorter than QBPCI for moderate sample size $n$; see Tables 1–4.

### Interval censored Case 1 data

In Banerjee and Wellner (2005), simulation results on the empirical likelihood-based confidence intervals (abbreviated as BWCI) for the median with interval censored Case 1 data (1.4) are presented in their Table 4, where $X$ and $Y$ both have $\text{Exp}(1)$ distribution. Here, we include some of our simulation results on smoothed WELRCI and non-smoothed WELRCI0 (see Remark 2) in Table 5 and compare them with BWCI, noting that for WELRCI0, high-order expansion is not relevant since smoothing of $\hat{F}_n$ is not used; see Remark 4. In Table 5, the simulation results for BWCI are taken directly from Table 4 of Banerjee and Wellner (2005). In our simulation studies, $c_n$ in (3.4) is determined according to the procedure described in Section 4, where for $k$-WELRCI or 1-WELRCI0, we use the following adaptive choice of the bootstrap sample size $n_b$ for the $m$ out of $n$ bootstrap to estimate $\rho_{n,\alpha}^{(k)}$:

$$n_b = \underset{\sqrt{n} \leq b_j \leq n^\gamma}{\arg\min} |\xi_{b_j,\alpha}^{(k)*} - \xi_{b_{j-1},\alpha}^{(k)*}|, \tag{5.1}$$

where $\gamma = 0.99, b_j = (b_{j-1} + d)$ with $b_0 = \sqrt{n}$ and $\xi_{b_j,\alpha}^{(k)*}$ is the $(1-\alpha)100$th percentile of $[b_j^{1/3}(\hat{\eta}_{b_j}^* - q)]^2/\hat{\mu}_2$. This selection method for the bootstrap sample sizes is a simple version of those in Götze and Račkauskas (2001) or Bickel and Sakov (2008) and our experience shows that this selected $n_b$ provides much more stable performance than simply using, say, $n_b = \sqrt{n}$ in simulation studies.

**Table 3.** 90% confidence intervals for $\theta_0 = F_0^{-1}(q)$ with exponential doubly censored data $X \sim \text{Exp}(1)$, $Y \sim \text{Exp}(3)$, $Z = (2/3)Y - 2.5$; percentage of $\delta$: $\delta = 1:56.0\%$; $\delta = 2: 24.9\%$; $\delta = 3: 19.1\%$

| | | $n=50$ | | $n=100$ | | $n=200$ | |
|---|---|---|---|---|---|---|---|
| $q$ & $\theta_0$ | Method | Coverage % | Average length (s.d.) | Coverage % | Average length (s.d.) | Coverage % | Average length (s.d.) |
| $q=0.250$ | 1-WELRCI | 89.5 | 0.306 (0.104) | 90.3 | 0.214 (0.063) | 89.9 | 0.153 (0.038) |
| $\theta_0 = 0.288$ | 2-WELRCI | 89.8 | 0.308 (0.105) | 90.6 | 0.215 (0.063) | 89.9 | 0.153 (0.038) |
| | SQBPCI | 89.4 | 0.298 (0.102) | 90.2 | 0.212 (0.062) | 89.9 | 0.152 (0.037) |
| | QBPCI | 89.6 | 0.324 (0.116) | 89.0 | 0.221 (0.068) | 90.2 | 0.155 (0.039) |
| $q=0.500$ | 1-WELRCI | 90.6 | 0.531 (0.175) | 90.3 | 0.375 (0.105) | 89.9 | 0.264 (0.063) |
| $\theta_0 = 0.693$ | 2-WELRCI | 90.0 | 0.520 (0.174) | 89.7 | 0.371 (0.104) | 89.8 | 0.263 (0.063) |
| | SQBPCI | 89.3 | 0.514 (0.169) | 89.4 | 0.369 (0.103) | 89.1 | 0.263 (0.063) |
| | QBPCI | 89.7 | 0.557 (0.195) | 91.6 | 0.389 (0.112) | 90.5 | 0.269 (0.065) |
| $q=0.750$ | 1-WELRCI | 89.8 | 1.017 (0.463) | 89.3 | 0.711 (0.264) | 90.1 | 0.497 (0.144) |
| $\theta_0 = 1.386$ | 2-WELRCI | 90.0 | 1.025 (0.466) | 89.4 | 0.714 (0.266) | 90.2 | 0.498 (0.144) |
| | SQBPCI | 88.7 | 0.982 (0.422) | 87.6 | 0.706 (0.263) | 88.7 | 0.496 (0.142) |
| | QBPCI | 90.3 | 1.155 (0.621) | 89.6 | 0.751 (0.294) | 91.4 | 0.516 (0.153) |

**Table 4.** 90% confidence intervals for $\theta_0 = F_0^{-1}(q)$ with chi-squared doubly censored data $X \sim \chi^2(1)$, $Y \sim \text{Exp}(3)$, $Z = (2/3)Y - 2.5$; percentage of $\delta$: $\delta = 1: 57.2\%$; $\delta = 2: 22.5\%$; $\delta = 3: 20.3\%$

| | | $n=50$ | | $n=100$ | | $n=200$ | |
|---|---|---|---|---|---|---|---|
| $q$ & $\theta_0$ | Method | Coverage % | Average length (s.d.) | Coverage % | Average length (s.d.) | Coverage % | Average length (s.d.) |
| $q=0.250$ | 1-WELRCI | 89.5 | 0.195 (0.093) | 88.8 | 0.139 (0.055) | 88.3 | 0.098 (0.030) |
| $\theta_0 = 0.102$ | 2-WELRCI | 90.0 | 0.196 (0.093) | 89.0 | 0.139 (0.055) | 88.3 | 0.098 (0.030) |
| | SQBPCI | 89.1 | 0.190 (0.090) | 88.9 | 0.136 (0.053) | 88.0 | 0.097 (0.029) |
| | QBPCI | 88.5 | 0.218 (0.104) | 88.2 | 0.147 (0.057) | 87.8 | 0.100 (0.031) |
| $q=0.500$ | 1-WELRCI | 89.2 | 0.543 (0.224) | 88.4 | 0.389 (0.128) | 90.1 | 0.278 (0.074) |
| $\theta_0 = 0.455$ | 2-WELRCI | 88.0 | 0.533 (0.222) | 87.8 | 0.386 (0.128) | 89.9 | 0.277 (0.074) |
| | SQBPCI | 87.1 | 0.521 (0.219) | 87.8 | 0.381 (0.126) | 88.8 | 0.275 (0.073) |
| | QBPCI | 89.3 | 0.586 (0.249) | 88.1 | 0.406 (0.140) | 89.8 | 0.284 (0.078) |
| $q=0.750$ | 1-WELRCI | 85.4 | 1.368 (0.754) | 86.9 | 0.928 (0.371) | 89.2 | 0.674 (0.208) |
| $\theta_0 = 1.323$ | 2-WELRCI | 85.9 | 1.378 (0.754) | 87.1 | 0.931 (0.371) | 89.3 | 0.675 (0.208) |
| | SQBPCI | 85.3 | 1.325 (0.672) | 86.9 | 0.929 (0.374) | 88.3 | 0.672 (0.209) |
| | QBPCI | 88.2 | 1.575 (0.924) | 88.3 | 1.006 (0.425) | 88.8 | 0.707 (0.227) |

From Table 5, we see that for interval censored Case 1 data, smoothed WELRCI perform well with the selected bootstrap sample size $n_b$, by (5.1). Compared with BWCI, WELRCI have better coverage for moderate sample sizes (also better than WELRCI0), while the average length of the confidence intervals is not significantly greater. For large sample sizes, WELRCI and WELRCI0 have similar performances and they both have slight over-coverage, but the average length of these confidence intervals is quite a bit less than the average length of those of BWCI. However, it should be noted that for large

**Table 5.** 95% confidence intervals for $\theta_0 = F_0^{-1}(q)$ with exponential interval censored data Case 1 $X \sim \text{Exp}(1)$, $Y \sim \text{Exp}(1)$; Percentage of $\delta$: $\delta = 1$: 50.0%; $\delta = 0$: 50.0%

| Sample Size | $q$ & $\theta_0$ Method | $q = 0.25$, Coverage percentage | $\theta_0 = 0.288$ Average length | $q = 0.50$, Coverage percentage | $\theta_0 = 0.693$ Average length | $q = 0.75$, Coverage percentage | $\theta_0 = 1.386$ Average length |
|---|---|---|---|---|---|---|---|
| $n = 50$ | 1-WELRCI0 | 61.4 | 0.385 | 95.0 | 1.162 | 81.9 | 1.321 |
| $d = 5$ in (5.1) | 3-WELRCI | 87.7 | 0.710 | 95.1 | 1.083 | 84.5 | 1.431 |
| | BWCI | – | – | 93.8 | 0.962 | – | – |
| $n = 100$ | 1-WELRCI0 | 82.0 | 0.380 | 96.1 | 0.936 | 94.8 | 1.489 |
| $d = 10$ in (5.1) | 3-WELRCI | 95.1 | 0.535 | 97.1 | 0.909 | 96.3 | 1.498 |
| | BWCI | – | – | 93.5 | 0.788 | – | – |
| $n = 200$ | 1-WELRCI0 | 94.0 | 0.360 | 98.2 | 0.665 | 98.4 | 1.212 |
| $d = 20$ in (5.1) | 3-WELRCI | 96.8 | 0.401 | 97.8 | 0.653 | 97.4 | 1.199 |
| | BWCI | – | – | 94.1 | 0.626 | – | – |
| $n = 500$ | 1-WELRCI0 | 97.7 | 0.244 | 97.3 | 0.391 | 97.0 | 0.722 |
| $d = 50$ in (5.1) | 3-WELRCI | 98.3 | 0.252 | 97.4 | 0.391 | 97.3 | 0.730 |
| | BWCI | – | – | 95.7 | 0.470 | – | – |
| $n = 1000$ | 1-WELRCI0 | – | – | 97.0 | 0.258 | – | – |
| $d = 100$ in (5.1) | 3-WELRCI | – | – | 97.2 | 0.259 | – | – |
| | BWCI | – | – | 95.0 | 0.367 | – | – |

sample size $n$, computing WELRCI or WELRCI0 is extremely time-consuming due to the use of (5.1).

In regard to selection criterion (4.3) for $k$, simulation results show that 1-WELRCI and 2-WELRCI perform similarly for right censored data and doubly censored data, while 3-WELRCI generally perform noticeably better than 1-WELRCI for interval censored data. But, for brevity, Table 5 does not include the simulation results of 1-WELRCI.

Overall, all simulation results presented here support Theorem 1 in Section 3 and the general methodology for implementation described in Section 4. The use of our rather simple version of smoothed quantile estimate $\hat{\theta} = \tilde{F}_n^{-1}(q)$ of (3.6) in the proposed WEL-RCI procedure performs very well in all of our simulation studies.

## Concluding remarks

We have shown that the proposed smoothed WELRCI provides a consistent likelihood-based interval estimate for quantiles with various types of censored data, including some of those that have not been previously studied in the literature. Compared with existing methods, smoothed WELRCI perform favorably in all available simulation studies. Moreover, the theoretical coverage accuracy equation (3.13) or (3.16) for smoothed WELRCI leads to the actual coverage accuracy result for right censored data and sheds light on further studies of coverage accuracy; see Remark 5. Finally, we note that the methods developed in this article can easily be used to construct WELRCI for survival probabilities, M-statistic, trimmed mean, etc., with different types of censored data.

## 6. Proofs

**Proof of Theorem 1(i).** Since $\tilde{F}_{\boldsymbol{p}}$ in (3.2) is an increasing function on $[0, W_m]$ with range $[0, 1]$, we have

$$r(\theta_0) = \sup\left\{ \prod_{i=1}^m (p_i/\hat{p}_i)^{n\hat{p}_i} \,\Big|\, \sum_{i=1}^m p_i U_i = 0, p_i \geq 0, \sum_{i=1}^m p_i = 1 \right\}, \tag{6.1}$$

where $U_i = [H_i(\boldsymbol{W}, \theta_0) - q]$. Note that from (3.2) and (3.6), we have for $\eta_0 = q$,

$$\sum_{i=1}^m \hat{p}_i U_i = (\hat{\eta} - q) = (\hat{\eta} - \eta_0) \quad \text{and} \quad \sum_{i=1}^m \hat{p}_i \hat{U}_i = 0, \tag{6.2}$$

where $\hat{U}_i = [H_i(\boldsymbol{W}, \hat{\theta}) - q]$. To get an expression for $r(\theta_0)$ in (6.1), we let

$$H(\boldsymbol{p}, \lambda, \gamma)$$
$$= n \sum_{i=1}^m \hat{p}_i (\log p_i - \log \hat{p}_i) - \lambda n \sum_{i=1}^m p_i U_i + \gamma \left( 1 - \sum_{i=1}^m p_i \right),$$

then,
$$\frac{\partial H}{\partial p_i} = \left(\frac{n\hat{p}_i}{p_i} - \lambda n U_i - \gamma\right) = 0 \quad \Rightarrow \quad \tilde{p}_i = \frac{n\hat{p}_i}{\lambda n U_i + \gamma}$$

and, in turn, the constraints $\sum_{i=1}^m \tilde{p}_i U_i = 0$ and $\sum_{i=1}^m \tilde{p}_i = \sum_{i=1}^m \hat{p}_i = 1$ give

$$n = \sum_{i=1}^m n\hat{p}_i = \sum_{i=1}^m \tilde{p}_i(\lambda n U_i + \gamma) = \gamma \quad \Rightarrow \quad \tilde{p}_i = \hat{p}_i/(1 + \lambda U_i) \tag{6.3}$$

and $\lambda$ should satisfy

$$0 = \sum_{i=1}^m \tilde{p}_i U_i = \sum_{i=1}^m \frac{\hat{p}_i U_i}{1 + \lambda U_i} \equiv g(\lambda). \tag{6.4}$$

The desired solutions of (6.4) are in the interval $(-U_{(m)}^{-1}, -U_{(1)}^{-1})$ because, in (6.3), we require $0 < 1 + \lambda U_i, 1 \le i \le m$, and

$$r(\theta_0) \ge c_n \quad \Rightarrow \quad U_{(1)} < 0 < U_{(m)} \quad \Rightarrow \quad \hat{\mu}_2 > 0. \tag{6.5}$$

Since $g'(\lambda) < 0$, for $\lambda_0$ as the unique solution of $g(\lambda) = 0$ in the interval $(-U_{(m)}^{-1}, -U_{(1)}^{-1})$,

$$\log r(\theta_0) = -n \sum_{i=1}^m \hat{p}_i \log(1 + \lambda_0 U_i). \tag{6.6}$$

To study the asymptotic behavior of $\lambda_0$, we notice that from (6.2), (6.4) and (6.6), we have $g(\lambda_0) = 0, g(0) = \sum_{i=1}^m \hat{p}_i U_i = (\hat{\eta} - \eta_0)$ and

$$-(\hat{\eta} - \eta_0) = [g(\lambda_0) - g(0)] = g'(\xi)\lambda_0$$
$$= -\lambda_0 \sum_{i=1}^m \hat{p}_i U_i^2 (1 + \xi U_i)^{-2},$$

where $|\xi| \le |\lambda_0|$. Thus, noting that (3.2) implies

$$0 \le H_i(\boldsymbol{W}, x) \le 1 \quad \Rightarrow \quad \max_{1 \le i \le m} |U_i| = \max_{1 \le i \le m} |H_i(\boldsymbol{W}, \theta_0) - q| \le 1, \tag{6.7}$$

from $(1 + \xi U_i)^2 \le (1 + |\lambda_0|)^2$ we have

$$|\hat{\eta} - \eta_0| = |\lambda_0| \sum_{i=1}^m \hat{p}_i U_i^2 (1 + \xi U_i)^{-2} \ge |\lambda_0| \sum_{i=1}^m \hat{p}_i U_i^2 (1 + |\lambda_0|)^{-2}.$$

Since (6.7) implies $|\lambda_0| \le \max\{|U_{(1)}^{-1}|, |U_{(m)}^{-1}|\} = \max\{q, 1-q\} \equiv M_1$, then from (3.7), (AS3) and Theorem 4.2.2 of Chung (1974), we have that with probability 1,

$$|\lambda_0| \le \hat{\mu}_2^{-1}|\hat{\eta} - \eta_0|(1 + M_1)^2 \le \rho|\hat{\eta} - \eta_0|(1 + M_1)^2 \tag{6.8}$$

all but finitely often, where $\rho > 0$ is a constant. Hence, by (AS2),

$$\lambda_0 \xrightarrow{a.s.} 0 \qquad \text{as } n \to \infty. \tag{6.9}$$

To avoid overly tedious algebra, the rest of the proof will establish (3.8)–(3.10) for $k = 2$, while the method can easily be used for the case $k = 4$.

We let $h = g^{-1}$, then $g(\lambda_0) = 0$ and $g(0) = (\hat{\eta} - \eta_0)$ imply that $h(0) = \lambda_0$ and $h(\hat{\eta} - \eta_0) = 0$, respectively. Moreover, we have

$$\begin{aligned} h'(\hat{\eta} - \eta_0) &= \frac{1}{g'(0)} = -\frac{1}{\hat{\mu}_2}, \\ h''(\hat{\eta} - \eta_0) &= -\frac{g''(0)}{[g'(0)]^3} = \frac{2\hat{\mu}_3}{\hat{\mu}_2^3}, \\ h'''(\hat{\eta} - \eta_0) &= \frac{3[g''(0)]^2 - g'(0)g'''(0)}{[g'(0)]^5} = -\frac{6(2\hat{\mu}_3^2 - \hat{\mu}_2\hat{\mu}_4)}{\hat{\mu}_2^5}, \\ h^{(4)}(y) &= \frac{10g'(x)g''(x)g'''(x) - [g'(x)]^2 g^{(4)}(x) - 15[g''(x)]^3}{[g'(x)]^7}, \end{aligned} \tag{6.10}$$

where $x = h(y)$ and, from Taylor's expansion, we have

$$\begin{aligned} \lambda_0 = h(0) &= h(\hat{\eta} - \eta_0) - h'(\hat{\eta} - \eta_0)(\hat{\eta} - \eta_0) + \frac{1}{2} h''(\hat{\eta} - \eta_0)(\hat{\eta} - \eta_0)^2 \\ &\quad - \frac{1}{6} h'''(\hat{\eta} - \eta_0)(\hat{\eta} - \eta_0)^3 + \frac{1}{24} h^{(4)}(\xi)(\hat{\eta} - \eta_0)^4 \\ &= \frac{\hat{\eta} - \eta_0}{\hat{\mu}_2} + \frac{\hat{\mu}_3}{\hat{\mu}_2^3}(\hat{\eta} - \eta_0)^2 + \frac{2\hat{\mu}_3^2 - \hat{\mu}_2\hat{\mu}_4}{\hat{\mu}_2^5}(\hat{\eta} - \eta_0)^3 + R_h, \end{aligned} \tag{6.11}$$

where $R_h = \frac{1}{24} h^{(4)}(\xi)(\hat{\eta} - \eta_0)^4$, $|\xi| \leq |\hat{\eta} - \eta_0|$, satisfying $\zeta = h(\xi)$ with $|\zeta| \leq |\lambda_0|$. Since by (6.7) and (6.9) we have

$$\begin{aligned} \frac{\hat{\mu}_2}{(1+|\lambda_0|)^2} &\leq |g'(\zeta)| = \sum_{i=1}^{m} \frac{\hat{p}_i U_i^2}{(1+\zeta U_i)^2} \leq \frac{\hat{\mu}_2}{(1-|\lambda_0|)^2}, \\ |g''(\zeta)| &= 2 \left| \sum_{i=1}^{m} \frac{\hat{p}_i U_i^3}{(1+\zeta U_i)^3} \right| \leq \frac{2}{(1-|\lambda_0|)^3}, \\ |g'''(\zeta)| &= 6 \left| \sum_{i=1}^{m} \frac{\hat{p}_i U_i^4}{(1+\zeta U_i)^4} \right| \leq \frac{6}{(1-|\lambda_0|)^4}, \\ |g^{(4)}(\zeta)| &= 24 \left| \sum_{i=1}^{m} \frac{\hat{p}_i U_i^5}{(1+\zeta U_i)^5} \right| \leq \frac{24}{(1-|\lambda_0|)^5}, \end{aligned} \tag{6.12}$$

from (6.9)–(6.10) and (AS3), there exists a constant $M_h$ such that with probability 1,

$$|R_h| \leq M_h|\hat{\eta} - \eta_0|^4 \tag{6.13}$$

all but finitely often. Thus, from Taylor's expansion, we have in (6.6)

$$\begin{aligned}
-2\log r(\theta_0) &= 2n\sum_{i=1}^{m} \hat{p}_i \log(1+\lambda_0 U_i) \\
&= 2n\sum_{i=1}^{m} \hat{p}_i\{\lambda_0 U_i - \tfrac{1}{2}(\lambda_0 U_i)^2 + \tfrac{1}{3}(\lambda_0 U_i)^3 \\
&\qquad - \tfrac{1}{4}(\lambda_0 U_i)^4 + \tfrac{1}{5}(1+\zeta_i)^{-5}(\lambda_0 U_i)^5\} \\
&= 2n(\lambda_0(\hat{\eta}-\eta_0) - \tfrac{1}{2}\lambda_0^2\hat{\mu}_2 + \tfrac{1}{3}\lambda_0^3\hat{\mu}_3 - \tfrac{1}{4}\lambda_0^4\hat{\mu}_4) + R_1,
\end{aligned} \tag{6.14}$$

where $|\zeta_i| \leq |\lambda_0 U_i| \leq |\lambda_0|$ and $R_1 = \tfrac{2}{5}n\sum_{i=1}^{m}\hat{p}_i(1+\zeta_i)^{-5}(\lambda_0 U_i)^5$. Easily, from (6.7)–(6.8), we know that with probability 1,

$$|R_1| \leq \tfrac{2}{5}n(1-|\lambda_0|)^{-5}|\lambda_0|^5 \leq n|\hat{\eta}-\eta_0|^5 M_{R_1} \tag{6.15}$$

all but finitely often, where $0 < M_{R_1} < \infty$ is a constant.

By plugging in (6.11) and using tedious algebra to combine the terms with the same rate of convergence, we obtain in (6.14)

$$-2\log r(\theta_0) = \frac{n(\hat{\eta}-\eta_0)^2}{\hat{\mu}_2}\left\{1 + \frac{2\hat{\mu}_3}{3\hat{\mu}_2^2}(\hat{\eta}-\eta_0) + \frac{2\hat{\mu}_3^2 - \hat{\mu}_2\hat{\mu}_4}{2\hat{\mu}_2^4}(\hat{\eta}-\eta_0)^2\right\} + R_1 + R_2, \tag{6.16}$$

where $R_2$ represents all terms in (6.16) with order $n|\hat{\eta}-\eta_0|^j$, $j \geq 5$. Thus, from (6.7) and (6.13), we have that there exists a constant $M_{R_2}$ such that with probability 1,

$$|R_2| \leq n|\hat{\eta}-\eta_0|^5 M_{R_2} \tag{6.17}$$

all but for finitely often. Therefore, (3.8)–(3.10) follow from (6.14)–(6.17). □

**Proof of Theorem 1(ii).** From (3.3), (3.8) and (6.5), we have

$$\begin{aligned}
&P\{X_L \leq \theta_0 \leq X_U\} \\
&= P\{-2\log r(\theta_0) \leq -2\log c_n, \hat{\mu}_2 > 0\} \\
&= [P\{(A_n^{(k)} + (C_n)^2(\hat{\eta}-\eta_0)^{k+3}r_{n,k}) \leq \tilde{c}_n\} - P\{A_n^{(k)} \leq \tilde{c}_n\}] + P\{A_n^{(k)} \leq \tilde{c}_n\},
\end{aligned} \tag{6.18}$$

where for $\hat{U} = C_n(\hat{\eta}-\eta_0)/\sqrt{\hat{\mu}_2}$ and $l_k(\hat{U},\hat{\eta}) = \hat{U}^2\sum_{i=1}^{k}\hat{a}_i(\hat{\eta}-\eta_0)^i$, we have

$$A_n^{(k)} = n^{-1}(C_n)^2 B_n^{(k)} = \hat{U}^2 + l_k(\hat{U},\hat{\eta}), \qquad k=1,2,3,4. \tag{6.19}$$

Note that with $l_0(\hat{U}, \hat{\eta}) = 0$, from (6.7), (3.8)–(3.9) and (AS2)–(AS3), we have that with probability 1, $[\hat{U}^2 + l_k(\hat{U}, \hat{\eta})] \geq \frac{1}{2}\hat{U}^2$ all but finitely often. Thus, for $\hat{r}_{n,k} = \hat{\mu}_2 r_{n,k}$, from (3.8), we know that $[\hat{U}^2 + l_k(\hat{U}, \hat{\eta})] \leq (\tilde{c}_n + |\hat{\eta} - \eta_0|^{k+1}|\hat{r}_{n,k}|\hat{U}^2)$ implies

$$\tfrac{1}{2}\hat{U}^2 \leq [\hat{U}^2 + l_k(\hat{U}, \hat{\eta})] \leq (\tilde{c}_n + |\hat{\eta} - \eta_0|^{k+1}|\hat{r}_{n,k}|\hat{U}^2) \leq (\tilde{c}_n + |\hat{\eta} - \eta_0|^{k+1}\hat{U}^2 M_{r,k}),$$

which, by (AS2), implies that with probability 1, $\tfrac{1}{4}\hat{U}^2 \leq (\tfrac{1}{2} - |\hat{\eta} - \eta_0|^{k+1} M_{r,k})\hat{U}^2 \leq \tilde{c}_n$ all but finitely often. Letting $M_{G_0}$ be the upper bound of $G_0'$ on $[a,b]$, since $\hat{U}^2 \leq 4\tilde{c}_n$ implies

$$|\hat{U}| \leq 2\sqrt{\tilde{c}_n} \quad \text{and} \quad |C_n(\hat{\eta} - \eta_0)| \leq 2\sqrt{\tilde{c}_n}, \tag{6.20}$$

(3.11) follows from (6.18)–(6.20), $\tilde{c}_n = O_p(1)$, (AS4) and

$$\begin{aligned}
&|P\{(A_n^{(k)} + (C_n)^2(\hat{\eta} - \eta_0)^{k+3} r_{n,k}) \leq \tilde{c}_n\} - P\{A_n^{(k)} \leq \tilde{c}_n\}| \\
&= |P\{(\hat{U}^2 + l_k(\hat{U}, \hat{\eta}) + \hat{U}^2(\hat{\eta} - \eta_0)^{k+1}\hat{r}_{n,k}) \leq \tilde{c}_n\} \\
&\quad - P\{[\hat{U}^2 + l_k(\hat{U}, \hat{\eta})] \leq \tilde{c}_n\}| \\
&\leq P\{(\tilde{c}_n - |\hat{\eta} - \eta_0|^{k+1}|\hat{r}_{n,k}|\hat{U}^2) \leq [\hat{U}^2 + l_k(\hat{U}, \hat{\eta})] \\
&\quad \leq (\tilde{c}_n + |\hat{\eta} - \eta_0|^{k+1}|\hat{r}_{n,k}|\hat{U}^2)\} \\
&\leq P\{(\tilde{c}_n - 4\tilde{c}_n M_{r,k}|\hat{\eta} - \eta_0|^{k+1}) \leq A_n^{(k)} \leq (\tilde{c}_n + 4\tilde{c}_n M_{r,k}|\hat{\eta} - \eta_0|^{k+1})\} \\
&\leq F_{n,k}(\tilde{c}_n + 4\tilde{c}_n M_{r,k}(2\sqrt{\tilde{c}_n})^{k+1}(C_n)^{-(k+1)}) \\
&\quad - F_{n,k}(\tilde{c}_n - 4\tilde{c}_n M_{r,k}(2\sqrt{\tilde{c}_n})^{k+1}(C_n)^{-(k+1)}) \\
&\leq 2\|F_{n,k} - G_0\|_{[a,b]} + M_{G_0}(8\tilde{c}_n M_{r,k}(2\sqrt{\tilde{c}_n})^{k+1})(C_n)^{-(k+1)}. \quad \square
\end{aligned} \tag{6.21}$$

**Proof of Corollary 1.** Note that for right censored data, we have $C_n = \sqrt{n}$, and that from Ren (1997), Theorem 1, we have $|\hat{F}_n(\theta_0) - \tilde{F}_n(\theta_0)| \leq |\hat{F}_n(\theta_0) - \hat{F}_n(\theta_0-)| = O_{a.s.}(n^{-1})$, where $O_{a.s.}(1)$ is bounded with probability 1 for sufficiently large $n$. Thus,

$$A_n^{(0)} = B_n^{(0)} = \frac{n(\hat{\eta} - \eta_0)^2}{\hat{\mu}_2} = \hat{W}^2 + O_p(n^{-1/2}), \tag{6.22}$$

where $\eta_0 = q, \hat{\xi} = \hat{F}_n(\theta_0), \hat{\nu}_2 = \sum_{i=1}^m \hat{p}_i(I\{W_i \leq \theta_0\} - \eta_0)^2$ and $\hat{W} = \sqrt{n}(\hat{\xi} - \eta_0)/\sqrt{\hat{\nu}_2}$, because

$$\begin{aligned}
\hat{\nu}_2 &= (1 - 2\eta_0)\hat{\xi} + \eta_0^2 = (1 - 2q)\hat{F}_n(\theta_0) + q^2, \\
\hat{\mu}_2 &= q^2 - 2q\hat{\eta} + \hat{F}_n(\theta_0-) + |\hat{F}_n(\theta_0) - \hat{F}_n(\theta_0-)|O(1) = \hat{\nu}_2 + O_p(n^{-1}).
\end{aligned}$$

Via straightforward algebra, we have

$$\hat{W}^2 = \frac{n(\hat{\xi} - \eta_0)^2}{(1 - 2\eta_0)\hat{\xi} + \eta_0^2} = \gamma_n^2 \frac{(1 - \eta_0)}{\eta_0}\sigma^2 + \hat{W}^2(\hat{\xi} - \eta_0)\frac{\alpha_n}{(1 - \hat{\xi})^2} + \hat{W}^2 \beta_n, \tag{6.23}$$

where for $\hat{\sigma}_G, \hat{\sigma}$ and $\sigma$ as in Chen and Lo (1996),

$$\gamma_n = \frac{\sqrt{n}(\hat{\xi} - \eta_0)}{(1-\hat{\xi})\hat{\sigma}_G}, \qquad \Delta_1 = \hat{\sigma}^2 - \sigma^2, \qquad \Delta_2 = \frac{\hat{\sigma}^2}{\hat{\sigma}_G^2} - 1, \qquad \hat{\sigma}_G^2 = \frac{\hat{\sigma}^2}{1+\Delta_2},$$

$$\alpha_n = \frac{(1+\Delta_2)\sigma^2(\eta_0\hat{\xi} - 1 + \eta_0 - \eta_0^2)}{\eta_0(\sigma^2 + \Delta_1)}, \qquad \beta_n = \frac{\Delta_1 - \sigma^2\Delta_2}{\sigma^2 + \Delta_1}.$$

From (A.7) and equation (2.13) of Chen and Lo (1996), we have

$$P\{|\Delta_1| > n^{-1/2}(\log n)^{-1}\} = o(n^{-1/2}) \quad \text{and} \quad P\{|\Delta_2| > n^{-2/3}\} = o(n^{-1/2}).$$

Thus, there exists a constant $0 < M_\sigma < \infty$ such that $P\{|\alpha_n| > M_\sigma\} = o(n^{-1/2})$ and $P\{|\beta_n| > M_\sigma n^{-1/2}\} = o(n^{-1/2})$. Noting that (6.20) implies $|\hat{\eta}| \le (\eta_0 + 2\sqrt{\tilde{c}_n}n^{-1/2})$, from (6.22) and a similar argument used in (6.21) we have

$$|P\{(A_n^{(0)} + n(\hat{\eta} - \eta_0)^3 r_{n,0}) \le \tilde{c}_n\} - P\{A_n^{(0)} \le \tilde{c}_n\}| \le 2\|F_\gamma - G_Z\|_{[a,b]} + O(n^{-1/2}),$$

where for $d = \sigma\sqrt{(1-\eta_0)/\eta_0}$ and $Z$ as the standard normal random variable, $F_\gamma$ is the d.f. of $(d\gamma_n)^2$ and $G_Z$ is the d.f. of $(dZ)^2$. Hence, from (6.18), the proof follows by noting that Theorem 2 of Chen and Lo (1996) implies $\|F_\gamma - G_Z\|_{[a,b]} = O(n^{-1/2})$. □

# Appendix

**Proof of (3.3).** Since $\tilde{F}_{\boldsymbol{p}}$ in (3.2) is an increasing function on $[0, W_m]$, it is easy to show that $\tau(\boldsymbol{p}) = \tilde{F}_{\boldsymbol{p}}^{-1}(q)$ is continuous in $\boldsymbol{p}$. Also, note that $S_n$ can be expressed as $S_n = \{\tau(\boldsymbol{p})|\boldsymbol{p} \in E_n\}$, where $E_n = \{\boldsymbol{p}|p_i \ge 0, \sum_{i=1}^m p_i = 1, \prod_{i=1}^m (p_i/\hat{p}_i)^{n\hat{p}_i} \ge c_n\}$ is a compact and convex set in $\mathbb{R}^m$. Since $\tau$ is continuous, from Royden (1968), page 158, we know that $S_n$ is a compact set in $\mathbb{R}$. Since convexity implies connectivity, from Royden (1968), pages 152–153, we know that $S_n$ is either an interval or a single point. Since $S_n$ is compact, we know that $S_n$ must be a closed interval $[X_L, X_U]$ with $X_L$ and $X_U$ given by (3.4). Next, we show (3.3) by denoting $E_0 = \{\boldsymbol{p}|\tau(\boldsymbol{p}) = \theta_0, p_i \ge 0, \sum_{i=1}^m p_i = 1\}$.

Assume $r(\theta_0) \ge c_n$. Since $\tau$ is continuous and $\{\boldsymbol{p}|p_i \ge 0, \sum_{i=1}^m p_i = 1\}$ is a compact set, we know that $E_0$ is compact and is not empty by (3.1). Thus, (3.1) and (3.4) imply that $X_L \le \theta_0 \le X_U$.

Assume $X_L \le \theta_0 \le X_U$. Since $\tau$ is continuous with $X_L$ and $X_U$ as the lower and upper bound on $E_n$, respectively, we know that from the intermediate value theorem, there exists $\boldsymbol{p}_0 \in E_n$ such that $\tau(\boldsymbol{p}_0) = \theta_0$. Hence, (3.1) implies $r(\theta_0) \ge c_n$. □

# Acknowledgements

This research was supported in part by NSF Grants DMS-02-04182 and DMS-06-04488. The author thanks Peter Bickel for discussions and conversations while this manuscript was being prepared. The author also thanks the Editor, the Associate Editor and the referee for their valuable comments and suggestions on the earlier draft of this paper.